\documentclass{article}

\usepackage{amsmath}
\usepackage{amsfonts}
\usepackage{amssymb}
\usepackage{url}
\usepackage{color}
\usepackage{amssymb}
\usepackage{enumerate}
\usepackage{graphics}
\usepackage{graphicx}

\newtheorem{theorem}{Theorem}
\newtheorem{lemma}[theorem]{Lemma}

\newtheorem{corollary}[theorem]{Corollary}

\newcommand{\Q}{{\mathbb Q}}
\newcommand{\Z}{{\mathbb Z}}

\newcommand{\R}{\mathbb R}

\def\Graver{{\cal G}}

\def\Orthant_j{{\mathcal O}_{j}}

\def\veb{{\ve b}}

\def\ved{{\ve d}}
\def\veg{{\ve g}}
\def\ves{{\ve s}}
\def\veu{{\ve u}}
\def\vev{{\ve v}}
\def\vex{{\ve x}}
\def\vey{{\ve y}}
\def\vez{{\ve z}}

\usepackage{enumerate}
\def\ve#1{\mathchoice{\mbox{\boldmath$\displaystyle\bf#1$}}
{\mbox{\boldmath$\textstyle\bf#1$}}
{\mbox{\boldmath$\scriptstyle\bf#1$}}
{\mbox{\boldmath$\scriptscriptstyle\bf#1$}}}

\newcommand{\boproof}{\textbf{Proof.} }
\newcommand{\eoproof}{\hspace*{\fill} $\square$ \vspace{5pt}}

\begin{document}
\setlength{\parindent}{0pt} \setlength{\parskip}{2ex plus 0.4ex
minus 0.4ex}

\title{Nash-equilibria and $N$-fold integer programming}
\author{{\bf Raymond~Hemmecke}\\University of Magdeburg, Germany\\
\\
{\bf Shmuel Onn}\\Technion, Haifa, Israel\\
\\
{\bf Robert Weismantel}\\University of Magdeburg, Germany}
\date{}

\maketitle

\begin{abstract}
Inspired by a paper of R.~W.~Rosenthal, we investigate generalized Nash-equilibria of integer programming games. We show that generalized Nash-equilibria always exist and are related to an optimal solution of a so-called $N$-fold integer program. This link allows us to establish some polynomial time complexity results about solving this optimization problem and its inverse counter-part.
\end{abstract}

\section{A class of IP games}

There are $N$ players. For $i=1,\ldots,N$, let $A^i\in\Z^{d^i\times n}$, $\veb^i\in\Z^{d^i}$, and $\veu^i\in\Z^n_+$. Each player $i$ selects a solution $\vex^i$ to the linear Diophantine system $A^i\vex=\veb^i, 0\leq\vex\leq\veu^i, \vex\in\Z^n$, which is referred to as a pure strategy of player $i$. The costs to $i$ for playing his strategy depend on the strategies chosen by the other players. A player is called \emph{satisfied} if his selected strategy is cost-minimal given the strategies of all other players. A \emph{pure Nash-equilibrium} is an assignment of strategies to players such that all players are satisfied simultaneously. Note that in the described situation each player can choose his strategy independently of the strategies of the other players. In contrast to this, one could impose a coupling constraint on the possible strategies for all players together, say $\sum_{i=1}^N B^i\vex^i\leq\veb^0$ for matrices $B^i\in\Z^{m\times n}$, $i=1,\ldots,N$, and for some $\veb^0\in\Z^m$. Again, a player is called \emph{satisfied} if his selected strategy is cost-minimal given the strategies of all other players. Clearly, his strategy has to obey the coupling constraint resulting from fixing the strategies of the other players. A \emph{generalized Nash-equilibrium} is an assignment of strategies to players such that all players are satisfied simultaneously in this more general situation. Clearly, we can model the unrestricted situtation by simply putting all $B^i$ to zero matrices.

Let us be more specific now. Let there be given functions $c_j:\R_+\rightarrow\R_+$ that evaluate the cost of element $j=1,\ldots,n$. In our setting, we have costs $c_j(y)$ if all the players together use $y=\sum_{i=1}^N x^i_j$ units of element $j$. In formulae, the costs assigned to player $k$ are
\[
\sum_{j=1}^n c_j\left(\sum_{i=1}^N y^i_j\right)-\sum_{j=1}^n
c_j\left(\sum_{i=1,i\neq k}^N y^i_j\right)=\sum_{j=1}^n
c_j\left(y^k_j+\sum_{i=1,i\neq k}^N y^i_j\right)-\sum_{j=1}^n
c_j\left(\sum_{i=1,i\neq k}^N y^i_j\right).
\]
From now on, we assume that $c_j$ is a monotonously increasing function for all $j=1,\ldots,n$. The following result is an extension of a theorem by Rosenthal \cite{Rosenthal}.

\begin{theorem}
All games of the above class with at least one choice $(\vex^1,\ldots,\vex^N)$ of strategies that fulfill the coupling constraint $\sum_{i=1}^N B^i\vex^i\leq\veb^0$ possess at least one generalized Nash-equilibrium.
\end{theorem}

\boproof Consider the problem
\begin{equation}\label{Equation: Nash-IP}
\min\left\{\sum_{j=1}^n c_j\left(\sum_{i=1}^N x^i_j\right):
A^i\vex^i=\veb^i, 0\leq\vex^i\leq\veu^i, \sum_{i=1}^N B^i\vex^i\leq\veb^0, \vex^i\in\Z^n_+\right\}.
\end{equation}
By assumption, this problem has at least one solution. Moreover, as for $j=1,\ldots,n$, the functions $c_j$ are monotonously increasing, they are all bounded from below, and thus, Problem (\ref{Equation: Nash-IP}) has an optimal solution. We claim that any optimal solution to Problem (\ref{Equation: Nash-IP}) is a generalized Nash-equilibrium.

Let $(\vey^1,\ldots,\vey^N)$ be an optimal solution to Problem (\ref{Equation: Nash-IP}). For the purpose of deriving a contradiction, suppose that $(\vey^1,\ldots,\vey^N)$ is not a generalized Nash-equilibrium. Then there exists some player $i$, w.l.o.g. $i=1$, and a strategy $\vez^1$ such that $\vez^1$ has lower costs
\[
\sum_{j=1}^n c_j\left(z^1_j+\sum_{i=2}^N y^i_j\right)-\sum_{j=1}^n
c_j\left(\sum_{i=2}^N y^i_j\right)
\]
than the costs
\[
\sum_{j=1}^n c_j\left(y^1_j+\sum_{i=2}^N y^i_j\right)-\sum_{j=1}^n
c_j\left(\sum_{i=2}^N y^i_j\right)
\]
of strategy $\vey^1$ given fixed strategies $\vey^2,\ldots,\vey^N$. We obtain
\[
\sum_{j=1}^n c_j\left(z^1_j+\sum_{i=2}^N y^i_j\right)-\sum_{j=1}^n
c_j\left(\sum_{i=2}^N y^i_j\right)<\sum_{j=1}^n
c_j\left(y^1_j+\sum_{i=2}^N y^i_j\right)-\sum_{j=1}^n
c_j\left(\sum_{i=2}^N y^i_j\right),
\]
that is,
\[
\sum_{j=1}^n c_j\left(z^1_j+\sum_{i=2}^N y^i_j\right)<\sum_{j=1}^n
c_j\left(y^1_j+\sum_{i=2}^N y^i_j\right).
\]
Thus, $(\vez^1,\vey^2,\ldots,\vey^N)$ is a feasible solution to
Problem (\ref{Equation: Nash-IP}) with strictly smaller costs than
the solution $(\vey^1,\vey^2,\ldots,\vey^N)$, contracting
optimality. \eoproof

The following remark seems in order. Notice that the objective function in Problem(\ref{Equation: Nash-IP}) can be interpreted as the costs of a provider to account for realizing all the strategies. The costs assigned to a player, however, are marginal costs. If all $c_j$, $j=1,\ldots,n$, are nonlinear convex functions, the sum of costs for all players exceeds the costs for the provider, i.e., it generates a gain for the provider. The more players play the bigger this gain gets. This fact motivates the provider to have all his players being satisfied. Interestingly, each player is satisfied when they cooperate to minimize the costs of the provider.

\section{Complexity results for convex cost functions}

Let us assume in this section that all matrices $A^i$ are equal to $A$, that all matrices $B^i$ are equal to $B$, and that the $c_j:\R_+\rightarrow\R_+$, $j=1,\ldots,n$, are univariate \emph{convex} functions. This latter condition is quite natural, since it simply means that the difference functions $c_j(y+1)-c_j(y)$ are monotonously increasing, that is, the cost for one extra unit of component $j$ is monotonously increasing the more of component $j$ is used.

With the simplifications, Problem (\ref{Equation: Nash-IP}) has a nice combinatorial structure:
{\tiny
\begin{equation}\label{Equation: N-fold Nash}
\min\left\{\sum_{j=1}^n c_j(y_j): \left(
  \begin{array}{cccccc}
    I_n & I_n & \cdots & I_n & -I_n & 0 \\
    B   & B   & \cdots & B   &    0 & I_n     \\
    A & 0 & \cdots & 0         & 0 & 0 \\
    0 & A & \cdots & 0         & 0 & \\
    \vdots & \vdots & \ddots & \vdots & \vdots & \vdots\\
    0  & 0 & \cdots & A   & 0 & 0\\
  \end{array}
\right)\left(
\begin{array}{c}
\vex^1\\
\vex^2\\
\vex^3\\
\vdots\\
\vex^N\\
\vey\\
\ves\\
\end{array}
\right)=\left(
\begin{array}{c}
\ve 0\\
\veb^0\\
\veb^1\\
\veb^2\\
\vdots\\
\veb^N\\
\end{array}
\right),\left(
\begin{array}{c}
\vex^1\\
\vex^2\\
\vdots\\
\vex^N\\
\end{array}
\right)\leq\left(
\begin{array}{c}
\veu^1\\
\veu^2\\
\vdots\\
\veu^N\\
\end{array}
\right),\vex^1,\ldots,\vex^N,\vey\in\Z^n_+\right\}.
\end{equation}
}

This structure allows us to prove a complexity result about the computation of the generalized Nash-equilibrium. However, in order to prove this, we need to introduce the notion of the Graver basis of a matrix together with a few important facts about these bases and about $N$-fold integer programs.

For $D\in\Z^{d\times n}$ let $\Graver(D)$ denote the Graver basis of $D$. This Graver basis consists exactly of those nonzero elements $\vez$ in $\ker(A)\cap\Z^n=\{\vex:A\vex=\ve 0\}\cap\Z^n$ that cannot be written as a non-trivial sum $\vez=\veu+\vev$ where $\veu,\vev\in\ker(A)\cap\Z^n$ and where $\veu$ and $\vev$ have the same sign-pattern, that is, $u_jv_j\geq 0$ for all $j=1,\ldots,n$. For any matrix $D$, the Graver basis $\Graver(D)$ is a finite set \cite{Graver:75}, and it has the nice property that it provides an optimality certificate for the family of separable convex integer minimization problems $\min\left\{\sum_{j=1}^n f_j(x_j):D\vex=\ved, 0\leq\vex\leq\veu, \vex\in\Z^n_+\right\}$ for any choice of convex functions $f_j:\R\rightarrow\R$, $j=1,\ldots,n$, and for any choice of $\ved\in\Z^d$ \cite{Graver:75,Murota+Saito+Weismantel}. That is, for any choice of the input data and for any feasible solution $\vex^0$ of the problem, $\vex^0$ is optimal if and only if there is no $\veg\in\Graver(D)$ such that $\vex^0+\veg$ is feasible and $\sum_{i=j}^n f_j(x^0_j+g_j)<\sum_{i=j}^n f_j(x^0_j)$.

Having a finite optimality certificate available, one can augment any given feasible solution of the problem to an optimal solution by iteratively walking along directions provided by the optimality certificate. If a simple greedy augmentation is followed, only a polynomial number (in the encoding length of the input data) of Graver basis augmentation steps is needed to construct an initial feasible solution and to augment it to optimality \cite{DeLoera+Hemmecke+Onn+Weismantel,Hemmecke+Onn+Weismantel:oracle}. Thus, if we had a Graver basis $\Graver(D)$ that is of polynomial size in $\langle D\rangle$, this augmentation approach provides a polynomial time algorithm to solve $\min\left\{\sum_{j=1}^n f_j(x_j):D\vex=\ved, 0\leq\vex\leq\veu, \vex\in\Z^n_+\right\}$.

Using results of \cite{Hosten+Sullivant,Santos+Sturmfels}, it was shown in \cite{DeLoera+Hemmecke+Onn+Weismantel} that for fixed matrices $A$ and $B$ (of appropriate dimensions) the sizes and the encoding lengths of the Graver bases of the family of matrices
\[
  [A,B]^{(N)}:=
  \left(
    \begin{array}{cccccccc}
      B   &   B & \cdots & B \\
      A   &     &        &   \\
          &   A &        &   \\
          &     & \ddots &   \\
          &     &        & A \\
   \end{array}
 \right)
\]
increase only polynomially in $N$. Thus, if we choose $D=[A,B]^{(N)}$ the corresponding separable convex minimization problem, a so-called $N$-fold integer program, can be solved in polynomial time for fixed matrices $A$ and $B$ \cite{DeLoera+Hemmecke+Onn+Weismantel,Hemmecke+Onn+Weismantel:oracle}. The following lemma shows that for fixed matrices $A$ and $B$ also the Graver bases of Problem (\ref{Equation: N-fold Nash}) increase only polynomially in $N$. And thus Problem (\ref{Equation: N-fold Nash}) can be solved in polynomial time as claimed.

\begin{lemma}\label{Lemma: Problem matrix has poly-size Graver basis}
Let $A\in\Z^{d\times n}$ and $B\in\Z^{m\times n}$ be fixed matrices, and let
\[
  \overline{[A,B]}^{(N)}:=\left(
    \begin{array}{cccccc}
      I_n & I_n & \cdots & I_n & -I_n & 0 \\
      B   & B   & \cdots & B   &    0 & I_n     \\
      A & 0 & \cdots & 0         & 0 & 0 \\
      0 & A & \cdots & 0         & 0 & \\
      \vdots & \vdots & \ddots & \vdots & \vdots & \vdots\\
      0  & 0 & \cdots & A   & 0 & 0\\
    \end{array}
  \right).
\]
Then the sizes of the Graver bases of $\overline{[A,B]}^{(N)}$, $N=1,2,\ldots$, and their encoding lengths increase only polynomially in $N$.
\end{lemma}

\boproof
Consider the matrix
\[
C^{(N)}:=\left(
\begin{array}{cccccccccc}
  I_n & -I_n & 0   & I_n & -I_n & 0   & \cdots & I_n & -I_n & 0   \\
  0   &    0 & I_n & 0   &    0 & I_n & \cdots & 0   & 0    & I_n \\
  B   &    0 & 0   & B   &    0 & 0   & \cdots & B   & 0    & 0   \\
  A   &    0 & 0   & 0   &    0 & 0   & \cdots & 0   & 0    & 0   \\
  0   &    0 & 0   & A   &    0 & 0   & \cdots & 0   & 0    & 0   \\
  \vdots & \vdots & \vdots & \vdots & \vdots & \vdots & \ddots & \vdots & \vdots & \vdots \\
  0   &    0 & 0   & 0   &    0 & 0   & \cdots & A   & 0    & 0   \\
\end{array}
\right).
\]
As this is the matrix of an $N$-fold IP composed of block matrices $\left(\begin{smallmatrix}A & 0 & 0 \\\end{smallmatrix}\right)$ and $\left(\begin{smallmatrix}I_n & -I_n& 0\\0&0&I_n\\B&0&0\\\end{smallmatrix}\right)$, the sizes of the Graver bases of $C^{(N)}$, $N=1,2,\ldots$, and their encoding lengths increase only polynomially in $N$, given that $A$, $B$, and $n$ are constant.

The matrix $\overline{[A,B]}^{(N)}$ is composed out of columns of $C^{(N)}$ and thus any Graver basis element of $\overline{[A,B]}^{(N)}$ filled with zeros for the additional columns of $C^{(N)}$ defines a Graver basis element of $C^{(N)}$. Therefore, also the sizes of the Graver bases of $\overline{[A,B]}^{(N)}$, $N=1,2,\ldots$, and their encoding lengths increase only polynomially in $N$, given that $A$, $B$, and $n$ are constant. \eoproof

\begin{theorem}\label{Theorem: Complexity of N-fold Nash-IP}
Assume that $A\in\Z^{d\times n}$ and $B\in\Z^{m\times n}$ are fixed and assume that all $c_j:\R_+\rightarrow\R_+$, $j=1,\ldots,n$, are univariate convex functions given by evaluation oracles. Moreover, let $M$ be a bound for $c_j\left(\sum_{i=1}^N x^i_j\right)$ for all $j$ and for all feasible selections of strategies for all players. Then Problem (\ref{Equation: N-fold Nash}) can be solved, and hence a generalized Nash-equilibrium can be computed, in time that is polynomial in $N$, in $\langle M, \veb^0,\veb^1,\ldots,\veb^N,\veu^1,\ldots,\veu^N\rangle$, and in the number of calls of evaluation oracles.
\end{theorem}

\boproof From Lemma \ref{Lemma: Problem matrix has poly-size Graver basis}, we know that the Graver basis of the problem matrix $\overline{[A,B]}^{(N)}$ has a polynomial-size encoding length in the encoding length of the input data. The result now follows from Theorem $7$ in \cite{Hemmecke+Onn+Weismantel:oracle}, which bounds the number of augmentations so as to reach an optimal solution of Problem
(\ref{Equation: N-fold Nash}). \eoproof

\section{Inverse optimization problem}

A (generalized) Nash-equilibrium satisfies all players simultaneously. It requires that cost functions $c_j$, $j=1,\ldots,n$, are given. Conversely, one may ask for finding cost functions $c_j$, $j=1,\ldots,n$, that turn a given set $\{\vex^1,\ldots,\vex^N\}$ of strategies into a Nash-equilibrium. Such inverse optimization problems are fundamental for a variety of applications, such as quality control in production systems, portfolio optimization, or production capacity planning, see \cite{Zhang} for details and references.

The problem in this generality is hopeless unless some structure on the cost functions $c_j$, $j=1,\ldots,n$, can be assumed. We study the situation that each $c_j$ has the form $c_j(y)=\lambda_j f_j(y)$ for some fixed given convex function $f_j:\R_+\rightarrow\R_+$ and some $\lambda_j\in\R_+$ that can be adjusted.

We can abstract from this setting slightly and define a general inverse integer optimization problem as follows.

\begin{lemma}
Consider the problem (IIOP): Given a polyhedron $P=\{\vex:D\vex=\ved,\ve 0\leq \vex\leq\veu\}\subseteq\R^n$, $\vex^*\in P\cap\Z^n$ and convex functions $f_j:\R\rightarrow\R$. Do there exist $\lambda_j\in\R_+$, $j=1,\ldots,n$, not all simultaneously zero, such that $\vex^*$ minimizes $f(\vex):=\sum_{j=1}^n \lambda_j f_j(x_j)$ over all $\vex\in P\cap\Z^n$.

\begin{itemize}\label{Theorem: IIOP}
  \item[(a)] (IIOP) is in coNP.
  \item[(b)] (IIOP) is solvable in time polynomial in $\langle\Graver(D)\rangle$, where $\Graver(D)$ denotes the Graver basis of $D$.
\end{itemize}
\end{lemma}

\boproof By Theorem $3$ in \cite{Hemmecke+Onn+Weismantel:oracle}, $\vex^*$ minimizes $f(\vex)=\sum_{j=1}^n \lambda_j f_j(x_j)$ over all $\vex\in P\cap\Z^n$ if and if only if $f(\vex^*+\veg)\geq f(\vex^*)$ for all $\veg\in\Graver(D)$ such that $\vex^*+\veg\in P$. Clearly, this can be checked in polynomial time in $\langle\Graver(D)\rangle$, showing (b).

To show (a), observe that this criterion is equivalent to asking for a solution to the following system in variables $\lambda_1,\ldots,\lambda_n$:
\begin{eqnarray*}
\sum_{j=1}^n \lambda_j f_j(x^*_j+g_j)& \geq & \sum_{j=1}^n \lambda_j
f_j(x^*_j) \;\;\;\text{ for all } \veg\in\Graver(D) \text{ such that
}
\vex^*+\veg\in P,\\
\lambda_1,\ldots,\lambda_n & \geq & 0,\\
\sum_{j=1}^n \lambda_j & > & 0,
\end{eqnarray*}
or equivalently, letting $H=\{\veg\in\Graver(D) \text{ such that }
\vex^*+\veg\in P\}$,
\begin{eqnarray*}
\sum_{j=1}^n \left[f_j(x^*_j+g_j)-f_j(x^*_j)\right] \lambda_j & \geq
& 0 \;\;\;\text{ for all } \veg\in H,\\
\lambda_1,\ldots,\lambda_n & \geq & 0,\\
\sum_{j=1}^n \lambda_j & > & 0.
\end{eqnarray*}
If this system has no solution, then there exists a subsystem
\begin{eqnarray*}
\sum_{j=1}^n \left[f_j(x^*_j+g_j)-f_j(x^*_j)\right] \lambda_j & \geq
& 0 \;\;\;\text{ for all } \veg\in H',\\
\lambda_1,\ldots,\lambda_n & \geq & 0,\\
\sum_{j=1}^n \lambda_j & > & 0,
\end{eqnarray*}
for $H'\subseteq H$ of at most $|H'|\leq n+1$ inequalities having no solution. (This is a consequence of Caratheodory's theorem and the Farkas-Lemma, see p. 94 in \cite{Schrijver:86}.)

Now we can apply standard linear programming technology: From the Farkas-Lemma, it follows that the latter system is not solvable if and only if there exists a vector $\vev=(v_{\veg})_{\veg\in H'}\in\Q_+^{|H'|}$ such that
\[
\sum_{\veg\in H'} \vev_{\veg}
\left[f_j(x^*_j+g_j)-f_j(x^*_j)\right]<0\;\;\;\text{ for all }
j=1,\ldots,n.
\]
This is a polynomially sized certificate for the ``No.''-answer to (IIOP), showing (a). \eoproof

Notice that Lemma \ref{Theorem: IIOP} (b) implies that if a Graver basis happens to be of polynomial size in the encoding length of the input data $\langle D,u\rangle$ of the integer program, then the corresponding optimization problem (IIOP) can be solved in polynomial time in $\langle D,u\rangle$.

Returning to Problem (\ref{Equation: N-fold Nash}), we want to solve its inverse counter-part. Note that again we assume that each $c_j$ has the form
$c_j(y)=\lambda_j f_j(y)$ for some fixed given convex function $f_j:\R_+\rightarrow\R_+$ and some $\lambda_j\in\R_+$ that can be adjusted.

\begin{corollary}
Let the matrices $A$ and $B$ be fixed and assume that all $f_j:\R_+\rightarrow\R_+$, $j=1,\ldots,n$, are univariate convex functions given by evaluation oracles. Moreover, let $M$ be a bound for $c_j\left(\sum_{i=1}^N x^i_j\right)$ for all $j$ and for all feasible selections of strategies for all players. Then the inverse integer optimization problem corresponding to Problem (\ref{Equation: N-fold Nash}) is solvable in time polynomial in $N$, in $\langle M, \vex^*, \veb^0,\veb^1,\ldots, \veb^N, \veu^1,\ldots, \veu^N\rangle$, and in the number of calls of evaluation oracles.
\end{corollary}

\boproof From Lemma \ref{Lemma: Problem matrix has poly-size Graver basis}, we know that the Graver basis of the problem matrix $\overline{[A,B]}^{(N)}$ of Problem (\ref{Equation: N-fold Nash}) is of polynomial size in $N$. The result thus follows from Theorem \ref{Theorem: IIOP}(b). \eoproof

\section{Extension of results}

Let us call the pair of matrices $(A^i,B^i)$ the \emph{type} of player $i$. While in Problem (\ref{Equation: Nash-IP}) all players could have a different type, we assumed in Problem (\ref{Equation: N-fold Nash}) that all players were of the same type $(A,B)$.
Let us assume now that there is a \emph{constant} number of possible player types $(A^i,B^i)$ represented by pairs of matrices $(\bar{A}^i,\bar{B}^i)$, $i=1,\ldots,t$, as in Problem (\ref{Equation: Nash-IP}). Despite the more general problem, the complexity results of the last two sections still hold.

\begin{theorem}
Assume that there is a \emph{constant} number $t$ of possible player types represented by pairs $(\bar{A}^i,\bar{B}^i)$, $i=1,\ldots,t$, of fixed matrices $\bar{A}^i$ and $\bar{B}^i$, $i=1,\ldots,t$. For $i=1,\ldots,N$, select one of these pairs $(\bar{A}^i,\bar{B}^i)$, $i=1,\ldots,t$, as the  type $(A^i,B^i)$ of player $i$. Finally, let $M$ be a bound for $c_j\left(\sum_{i=1}^N x^i_j\right)$ for all $j$ and for all feasible selections of strategies for all players.
\begin{itemize}
  \item[(a)] Problem (\ref{Equation: Nash-IP}) can be solved in time polynomial in $N$, in $\langle M,\veb^0,\veb^1,\ldots,\veb^N,\veu^1,\ldots,\veu^N\rangle$, and in the number of calls of evaluation oracles.
  \item[(b)] Assume that each $c_j$ has the form $c_j(y)=\lambda_j f_j(y)$ for some fixed convex function $f_j:\R_+\rightarrow\R_+$ given by evaluation oracles and some $\lambda_j\in\R_+$ that can be adjusted. Then the inverse integer optimization problem corresponding to Problem (\ref{Equation: Nash-IP}) is solvable in time polynomial in $N$, in $\langle M,\vex^*,\veb^0,\veb^1,\ldots,\veb^N,\veu^1,\ldots,\veu^N\rangle$, and in the number of calls of evaluation oracles.
\end{itemize}
\end{theorem}

\boproof
As in the proofs above, it suffices to show that under our assumptions the sizes of the Graver bases of the problem matrices of Problem (\ref{Equation: Nash-IP}) increase only polynomially with $N$. This can be seen by considering the matrices
\[
\left[\left(\begin{smallmatrix}
\bar{A}^1 & 0         & \ldots & 0\\
0         & \bar{A}^2 & \ldots & 0\\
0         & 0         & \ddots & 0\\
0         & 0         & \ldots & \bar{A}^t\\
\end{smallmatrix}
\right),\left(\begin{smallmatrix}
\bar{B}^1 & \bar{B}^2 & \ldots & \bar{B}^t\\
\end{smallmatrix}
\right)\right]^{(N)}=
\left(\begin{array}{ccccccccc}
\bar{B}^1 & \bar{B}^2 & \ldots & \bar{B}^t & \ldots & 0         & 0         & 0      & 0         \\
\bar{A}^1 & 0         & \ldots & 0         & \ldots & 0         & 0         & 0      & 0         \\
0         & \bar{A}^2 & \ldots & 0         & \ldots & 0         & 0         & 0      & 0         \\
0         & 0         & \ddots & 0         & \ldots & 0         & 0         & 0      & 0         \\
0         & 0         & \ldots & \bar{A}^t & \ldots & 0         & 0         & 0      & 0         \\
\vdots    & \vdots    & \vdots & \vdots    & \ddots & \vdots    & \vdots    & \vdots & \vdots    \\
0         & 0         & 0      & 0         & \ldots & \bar{B}^1 & \bar{B}^2 & \ldots & \bar{B}^t \\
0         & 0         & 0      & 0         & \ldots & \bar{A}^1 & 0         & \ldots & 0         \\
0         & 0         & 0      & 0         & \ldots & 0         & \bar{A}^2 & \ldots & 0         \\
0         & 0         & 0      & 0         & \ldots & 0         & 0         & \ddots & 0         \\
0         & 0         & 0      & 0         & \ldots & 0         & 0         & \ldots & \bar{A}^t \\
\end{array}
\right),
\]
which are matrices of $N$-fold integer programs. Thus, as these $N$-fold matrices are built from fixed blocks, the sizes and encoding lengths their Graver bases increase only polynomially with $N$. As the problem matrix of Problem (\ref{Equation: Nash-IP}) can be obtained from the matrix by deleting columns and then zero rows, the sizes and the encoding lengths of the Graver bases of the problem matrices of Problem (\ref{Equation: Nash-IP}) increase only polynomially with $N$, and thus, (a) and (b) both follow. \eoproof

\end{document}